\documentclass{pspum-l}

\newtheorem{theorem}{Theorem}[section]
\newtheorem{lemma}[theorem]{Lemma}
\newtheorem{prop}[theorem]{Proposition}
\newtheorem{conj}[theorem]{Conjecture}
\newtheorem{problem}[theorem]{Problem}

\theoremstyle{definition}
\theoremstyle{remark}

\numberwithin{equation}{section}

\usepackage[all]{xy}
\usepackage{graphicx} 
\usepackage{tikz}
\usetikzlibrary{patterns}
\usetikzlibrary{decorations.pathreplacing}
\usetikzlibrary{arrows}
\usetikzlibrary{decorations.markings}

\usepackage{amsmath}
\usepackage{amssymb}
\usepackage{enumerate}
\usepackage{mathdots}
\usepackage{mathtools}
\usepackage{hyperref}
\usepackage{amsthm}

\allowdisplaybreaks

\usepackage{ytableau}
\ytableausetup{boxsize=0.8em}

\newcommand{\dfn}[1]{\textcolor{blue}{\emph{#1}}}

\begin{document}

\title{Order polynomial product formulas and poset dynamics}

\author{Sam Hopkins}
\address{Department of Mathematics, Howard University, Washington, DC}
\email{samuelfhopkins@gmail.com}
\thanks{The author was supported in part by NSF grant \#1802920.}

\begin{abstract}
We survey all known examples of finite posets whose order polynomials have product formulas, and we propose the heuristic that these are the same posets with good dynamical behavior. Here the dynamics in question are the actions of promotion on the linear extensions of the poset and rowmotion on the $P$-partitions of the poset.
\end{abstract}

\maketitle

\section{Introduction} \label{sec:intro}

A \dfn{linear extension} of a poset is a linear (i.e., total) ordering of its elements extending the partial ordering. Perhaps the single most important numerical invariant associated to a finite\footnote{All posets we consider will be finite and we will drop this adjective from now on.}  poset $P$ is its number $e(P)$ of linear extensions. Computing the number of linear extensions of a poset is a $\#\mathrm{P}$-complete problem~\cite{brightwell1991counting}, which means that we cannot hope for a ``good'' formula for this number in general. Nevertheless, especially within the context of algebraic combinatorics, there is great interest in obtaining good formulas -- in the best cases, product formulas -- for the number of linear extensions of special families of posets. The most famous such product formula is the celebrated hook-length formula~\cite{frame1954hook} for the number of linear extensions of a poset of Young diagram shape. There are also hook-length formulas for shifted shapes, and for rooted forests; and more generally, the $d$-complete posets of Proctor~\cite{proctor2014dcomplete, kim2019hook, naruse2019skew}, which include Young diagram shapes, shifted shapes, and rooted forests, have hook-length formulas counting their linear extensions. 

Here we will concentrate on a finer combinatorial invariant of a poset $P$ than its number of linear extensions: namely, its order polynomial. A \dfn{$P$-partition of height~$m$} is a weakly order-preserving map $P\to \{0,1,\ldots,m\}$; the \dfn{order polynomial} $\Omega_P(m)$ of~$P$ counts the number of $P$-partitions of height $m$.\footnote{Traditionally, as in~\cite[Chapter~3]{stanley2012ec1}, $P$-partitions are order-reversing and the order polynomial is what is $\Omega_P(m-1)$ in our notation. Our changes are superficial but lead to a cleaner presentation.}  The order polynomial is a polynomial in $m$ of degree $\#P$, with leading coefficient equal to~$e(P)/\#P!$. Posets with product formulas for their order polynomials are much rarer than posets with product formulas enumerating their linear extensions.

\medskip

We will review all known examples of posets with product formulas for their order polynomials below. However, our goal is not just to survey families of posets with order polynomial product formulas, but also to advertise an apparently powerful heuristic which says that:
\begin{center} {\bf the posets with order polynomial product formulas are \\ the same as the posets with good dynamical behavior}.\end{center}

Let us explain what we mean by poset dynamics. There are several constructions from algebraic combinatorics of interesting invertible operators acting on objects associated to a poset (e.g., linear extensions or $P$-partitions). And while these operators are defined for any poset, they tend to have good behavior (e.g., a small, predictable order and regular orbit structure) only for a select few families of posets. Here we will focus on two such operators: \dfn{promotion} and \dfn{rowmotion}. 

Promotion is an invertible operator acting on the linear extensions of any poset. It was first defined and studied by Sch\"{u}tzenberger~\cite{schutzenberger1972promotion}, in conjunction with a related involutive operator called evacuation. Sch\"{u}tzenberger's initial motivation for studying these operators was the relation between evacuation and the RSK algorithm~\cite{schutzenberger1972schensted} (see~\cite[Chapter~7, \S A1.2]{stanley1999ec2} for a modern account of this relation).

Rowmotion is an invertible operator acting on the order ideals of any poset, which has been studied by a number of authors over several decades~\cite{brouwer1974period, deza1990loops, fonderflaass1993orbits, cameron1995orbits}, with a renewed interest especially in the last 10 or so years because of a surprising connection to the combinatorics of root systems~\cite{panyushev2009orbits, striker2012promotion, williams2013cataland}. Einstein and Propp~\cite{einstein2013combinatorial} gave a piecewise-linear extension of rowmotion to the entire order polytope $\mathcal{O}(P)$ of a poset $P$. By identifying the points in $\frac{1}{m}\mathbb{Z}^{P}\cap \mathcal{O}(P)$ with $P$-partitions of height $m$, we thus obtain an action of rowmotion on these $P$-partitions. It is this (piecewise-linear) action of rowmotion on $P$-partitions which we will mostly be concerned with.

Then, the heuristic we are proposing is more precisely that the following three properties of a poset $P$ are related and tend to occur simultaneously:
\begin{enumerate}[(1)]
\item \label{item:prod_form} the order polynomial $\Omega_P(m)$ has a product formula;
\item \label{item:pro} promotion acting on the linear extensions of $P$ has good behavior;
\item \label{item:row} rowmotion acting on the $P$-partitions of height $m$ has good behavior, for all $m$.
\end{enumerate}
These three properties are not perfectly correlated, in that we have examples of posets which satisfy some but not all of them (although we know of no counterexamples to \eqref{item:row} $\Rightarrow$ \eqref{item:pro} $\Rightarrow$ \eqref{item:prod_form}). Nevertheless, we do think this heuristic is powerful, and we remark that it is powerful ``in both directions'': that is, both for finding posets with good dynamical behavior (see, e.g.,~\cite{hopkins2020promotion}), and for finding posets with order polynomial product formulas (see, e.g., \cite{hopkins2020plane}).

\subsection{A more detailed account of the heuristic} \label{subsec:heuristic}

Let us give a more detailed account of what~\eqref{item:prod_form},~\eqref{item:pro}, and~\eqref{item:row} mean, and how they are related.

For $\Omega_P(m)$ to have a product formula, ideally all of its roots should be integers. There are also some interesting examples where the roots are \emph{half}-integers, so we will consider this acceptable as well. Furthermore, $P$ should come in a family, with certain numerical parameters attached to it, and we should be able to write $\Omega_P(m)$ in a simple way as a product of rational expressions involving these parameters. Since $e(P)=\#P! \cdot \lim_{m\to \infty} \frac{\Omega_P(m)}{m^{\#P}}$, whenever we have a product formula for $\Omega_P(m)$ we also have one for $e(P)$.

As it turns out, we also always seem to get very nice $q$-analogs when $P$ has an order polynomial product formula. Define
\[\Omega_P(m;q)  \coloneqq \prod_{\alpha} \frac{(1-q^{\kappa(m-\alpha)})}{(1-q^{-\kappa\alpha})},\]
 where the product is over all roots $\alpha$ of $\Omega_P(m)$, with multiplicity, and $\kappa$ is $1$ if these roots are all integers and $2$ if they are half-integers. Then, miraculously, in the examples we observe that $\Omega_P(m;q)$ is (for nonnegative integers $m$) a polynomial in~$q$ with nonnegative integer coefficients, which at $q=1$ is equal to $\Omega_P(m)$. Similarly, define
\[e(P;q)  \coloneqq \prod_{j=1}^{\#P}(1-q^{j\kappa}) \cdot \lim_{m\to \infty} \Omega_P(m;q) = \prod_{j=1}^{\#P}(1-q^{j\kappa}) \prod_{\alpha}\frac{1}{(1-q^{-\kappa\alpha})}.\] Again, $e(P;q)$ is miraculously a $q$-analog of $e(P)$ in the sense that it is a polynomial in $q$ with nonnegative integer coefficients which at $q=1$ is equal to $e(P)$. (See~\cite{stanton1990fake} for discussion of when $q$-expressions of these forms are polynomials with nonnegative integer coefficients.)

Let $\mathcal{L}(P)$ denote the linear extensions of $P$ and $\mathrm{Pro} \colon \mathcal{L}(P)\to \mathcal{L}(P)$ denote promotion. When we formally define promotion and go over the basics below, we will see why $\mathrm{Pro}^{\#P}$ is the ``right power'' of promotion to look at to find good behavior. For promotion to have good behavior, ideally $\mathrm{Pro}^{\#P}$ is the identity. There are also some interesting examples where $\mathrm{Pro}^{\#P}$ is a non-identity involutive poset automorphism: these are exactly the examples where the roots of $\Omega_P(m)$ are half-integers.

Moreover, whenever promotion of $\mathcal{L}(P)$ has good behavior, $e(P;q)$ is apparently a cyclic sieving polynomial for the action of promotion. We will review the cyclic sieving phenomenon later, but for now the important remark is that if a polynomial with an expression as a product of ratios of $q$-numbers is a cyclic sieving polynomial for some cyclic action, then this action has a very regular orbit structure: for instance, this means there is a product formula enumerating every symmetry class.

Let $\mathcal{PP}^m(P)$ denote the height $m$ $P$-partitions and $\mathrm{Row} \colon \mathcal{PP}^m(P)\to \mathcal{PP}^m(P)$ denote rowmotion. Rowmotion only ever has good behavior when $P$ is graded (i.e., all maximal chains of $P$ have the same length). For a graded poset $P$ we use $r(P)$ to denote the rank of $P$ (i.e., the length of a maximal chain); we will see below why $\mathrm{Row}^{r(P)+2}$ is the ``right power'' of rowmotion to look at. For rowmotion to have good behavior, ideally $\mathrm{Row}^{r(P)+2}$ is the identity; but there are also interesting examples where it is an involutive automorphism, and again this happens when the roots of~$\Omega_P(m)$ are half-integers. Moreover, when rowmotion of $\mathcal{PP}^m(P)$ has good behavior, $\Omega_P(m;q)$ is apparently a cyclic sieving polynomial for this action.

\subsection{Whence all this?} 

Where does this heuristic comes from, and why might these properties of a poset be related? The short answer is that these properties are indicative of some connection of the poset to \emph{algebra}, especially, the representation theory of Lie algebras, Lie groups, Weyl groups, etc. 

For instance, it often happens that $\mathcal{PP}^m(P)$ indexes a basis of an irreducible representation of a Lie algebra, in which case we can compute $\Omega_P(m)$ using the Weyl dimension formula. Similarly, the $q$-analog $\Omega_P(m;q)$ can be obtained via a $q$-Weyl dimension formula for the principal specialization of the corresponding character. 

Furthermore, the actions of promotion and rowmotion routinely have nice algebraic models as well. In the simplest cases, there are in fact diagrammatic models where the action is realized as rotation; but there are also examples where sophisticated tools from algebra like crystals and  canonical bases are required. Indeed, these models are part of what make promotion and rowmotion so fascinating. We will review all known models of promotion and rowmotion below.

For the families of posets which have a direct connection to algebra, it is desirable to prove results uniformly, that is, without relying on classification theorems. We should also note that for some posets (such as the ``chain of V's''), the relevant algebra has apparently yet to be uncovered.

\subsection{The open problems} 

As we will see in the subsequent sections, simply by carrying out the program in Section~\ref{subsec:heuristic} for the posets known to have an order polynomial product formula, we obtain many intriguing dynamical conjectures. But our heuristic also presents a few ``meta-problems'': 

\begin{problem} \label{prob:connection}
Find formal relations between the properties~\eqref{item:prod_form},~\eqref{item:pro}, and~\eqref{item:row}.
\end{problem}

\begin{problem}
Find more examples of posets satisfying~\eqref{item:prod_form},~\eqref{item:pro}, and~\eqref{item:row}.
\end{problem}

\begin{problem}
Find a unified algebraic explanation for the good (enumerative and dynamical) behavior of the families of posets considered here.
\end{problem}

Regarding Problem~\ref{prob:connection}, we should remark that while there are several papers which study connections between promotion and rowmotion (see, e.g.,~\cite{striker2012promotion, dilks2017resonance, dilks2019rowmotion, bernstein2020promotion, bernstein2022promotion}), we know of none which discusses a direct connection between promotion of~$\mathcal{L}(P)$ and rowmotion of~$\mathcal{PP}^m(P)$.

\subsection{Acknowledgments} 
I thank Ira Gessel, Soichi Okada, Rebecca Patrias, Robert Proctor, Victor Reiner, Martin Rubey, Jessica Striker, Bruce Westbury, and Nathan Williams for useful discussion. 

\section{The posets} \label{sec:posets}

In this section we introduce the posets which have order polynomial product formulas. We assume the reader is familiar with poset basics as laid out for instance in~\cite[Chapter~3]{stanley2012ec1}. All the properties of posets we are interested in decompose in a natural way over disjoint unions, so we will only consider connected posets. These properties also evidently translate directly from a poset $P$ to its dual $P^*$. In fact, these properties (at least conjecturally) translate from a poset to any other poset with an isomorphic comparability graph (see~\cite{hopkins2019minuscule}). Therefore, we will not separately list posets with isomorphic comparability graphs to the ones below.

\subsection{Shapes and shifted shapes}

Several of the families of posets which have order polynomial product formulas will be Young diagram shapes or shifted shapes.

We assume the reader is familiar with the basics concerning partitions, shapes, and so on. We view a shape as a poset on its boxes with the partial order where $u\leq v$ means box $u$ is weakly northwest of box $v$. The poset objects associated to shapes have different traditional names: e.g., a linear extension is a standard Young tableau, a $P$-partition is a plane partition, etc.

\begin{figure}
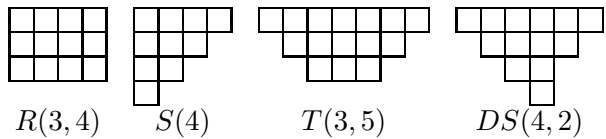

\begin{center}
\begin{tabular}{c c c c}
\ydiagram{4,4,4} &  \ydiagram{4,3,2,1} & \ydiagram{7,1+5,2+3} & \ydiagram{6,1+4,2+2,3+1} \\
$R(3,4)$ &$S(4)$ & $T(3,5)$ & $DS(4,2)$ \\
\end{tabular}
\end{center}
\caption{Examples of the families of shapes.} \label{fig:shapes}
\end{figure}

We now define the relevant families of shapes. The \dfn{rectangle} $R(a,b)$ is the shape for the partition $(b^a)$ using multiplicity notation.  The \dfn{staircase} $S(n)$ is the shape for the partition $(n,n-1,\ldots,1)$. The \dfn{shifted trapezoid} $T(a,b)$ is the shifted shape for the strict partition $(a+b-1,a+b-3,a+b-5,\ldots,b-a+1)$. The \dfn{shifted double staircase} $DS(n,k)$ is the shifted shape for the strict partition $(n,n-1,\ldots,1)+(k,k-1,\ldots,1)$. Observe that $T(n,n)=DS(n,n-1)$ and $T(n,n+1)=DS(n,n)$. Figure~\ref{fig:shapes} depicts examples of these shapes.

We also define the arithmetic progression $AP(M,d,\ell)$ to be the shape for the partition $(M-d,M-2d,\ldots,M-\ell d)$. Observe that $R(a,b)=AP(b,0,a)$ and $S(n)=AP(n+1,1,n)$. 

\subsection{Root posets}

The root posets are a very interesting family of posets coming from Lie theory, and some of them (specifically, the root posets of coincidental type) have order polynomial product formulas. We give only a very cursory account of root posets here; for a detailed account see~\cite{williams2013cataland} or~\cite[\S8]{hamaker2018doppelgangers}.

\begin{figure}
\begin{center}
\begin{tikzpicture}[scale=0.7]
	\node[shape=circle,fill=black,inner sep=1.5] (1) at (0,0) {};
	\node[shape=circle,fill=black,inner sep=1.5] (2) at (1,0) {};
	\node[shape=circle,fill=black,inner sep=1.5] (3) at (2,0) {};
	\node[shape=circle,fill=black,inner sep=1.5] (4) at (3,0) {};
	\node[shape=circle,fill=black,inner sep=1.5] (5) at (0.5,0.5) {};
	\node[shape=circle,fill=black,inner sep=1.5] (6) at (1.5,0.5) {};
	\node[shape=circle,fill=black,inner sep=1.5] (7) at (2.5,0.5) {};
	\node[shape=circle,fill=black,inner sep=1.5] (8) at (1,1) {};
	\node[shape=circle,fill=black,inner sep=1.5] (9) at (2,1) {};
	\node[shape=circle,fill=black,inner sep=1.5] (10) at (1.5,1.5) {};
	\draw [thick] (1) -- (5) -- (2) -- (6) -- (3) -- (7) -- (4);
	\draw [thick] (5) -- (8) -- (6) -- (9) -- (7);
	\draw [thick] (8) -- (10) -- (9);
	\draw [decorate,decoration={brace,amplitude=10pt,mirror},yshift=-2.5pt] (-0.1,0) -- (3.1,0) node [black,midway,yshift=-15pt]  {\footnotesize $n$};
	\node at (1.5,-1.5) {$\Phi^+(A_n)$};
\end{tikzpicture}  \vline  \begin{tikzpicture}[scale=0.7]
	\node[shape=circle,fill=black,inner sep=1.5] (1) at (0,0) {};
	\node[shape=circle,fill=black,inner sep=1.5] (2) at (1,0) {};
	\node[shape=circle,fill=black,inner sep=1.5] (3) at (2,0) {};
	\node[shape=circle,fill=black,inner sep=1.5] (4) at (3,0) {};
	\node[shape=circle,fill=black,inner sep=1.5] (5) at (0.5,0.5) {};
	\node[shape=circle,fill=black,inner sep=1.5] (6) at (1.5,0.5) {};
	\node[shape=circle,fill=black,inner sep=1.5] (7) at (2.5,0.5) {};
	\node[shape=circle,fill=black,inner sep=1.5] (8) at (1,1) {};
	\node[shape=circle,fill=black,inner sep=1.5] (9) at (2,1) {};
	\node[shape=circle,fill=black,inner sep=1.5] (10) at (1.5,1.5) {};
	\node[shape=circle,fill=black,inner sep=1.5] (11) at (0,1) {};
	\node[shape=circle,fill=black,inner sep=1.5] (12) at (0.5,1.5) {};
	\node[shape=circle,fill=black,inner sep=1.5] (13) at (0,2) {};
	\node[shape=circle,fill=black,inner sep=1.5] (14) at (1,2) {};
	\node[shape=circle,fill=black,inner sep=1.5] (15) at (0.5,2.5) {};
	\node[shape=circle,fill=black,inner sep=1.5] (16) at (0,3) {};
	\draw [thick] (1) -- (5) -- (2) -- (6) -- (3) -- (7) -- (4);
	\draw [thick] (5) -- (8) -- (6) -- (9) -- (7);
	\draw [thick] (8) -- (10) -- (9);
	\draw [thick] (5) -- (11) -- (12) -- (8);
	\draw [thick] (13) -- (12) -- (14) -- (10);
	\draw [thick] (13) -- (15) -- (14);
	\draw [thick] (15) -- (16);		
	\draw [decorate,decoration={brace,amplitude=10pt,mirror},yshift=-2.5pt] (-0.1,0) -- (3.1,0) node [black,midway,yshift=-15pt]  {\footnotesize $n$};
	\node at (1.5,-1.5) {$\Phi^+(B_n)\simeq \Phi^+(C_n)$};
\end{tikzpicture} \vline \;\; \begin{tikzpicture}[scale=0.7]
	\node[shape=circle,fill=black,inner sep=1.5] (1) at (0,0) {};
	\node[shape=circle,fill=black,inner sep=1.5] (2) at (1,0) {};
	\node[shape=circle,fill=black,inner sep=1.5] (3) at (2,0) {};
	\node[shape=circle,fill=black,inner sep=1.5] (4) at (3,0) {};
	\node[shape=circle,fill=black,inner sep=1.5] (5) at (0.5,0.5) {};
	\node[shape=circle,fill=black,inner sep=1.5] (6) at (1.5,0.5) {};
	\node[shape=circle,fill=black,inner sep=1.5] (7) at (2.5,0.5) {};
	\node[shape=circle,fill=black,inner sep=1.5] (8) at (1,1) {};
	\node[shape=circle,fill=black,inner sep=1.5] (9) at (2,1) {};
	\node[shape=circle,fill=black,inner sep=1.5] (11) at (0,1) {};
	\node[shape=circle,fill=black,inner sep=1.5] (12) at (0.5,1.5) {};
	\node[shape=circle,fill=black,inner sep=1.5] (13) at (0,2) {};
	\draw [thick] (1) -- (5) -- (8) -- (7);
	\draw [thick] (2) -- (6) -- (3) -- (7) -- (4);
	\draw [thick] (6) -- (9) -- (7);
	\draw [thick] (3) -- (5) -- (11) -- (6);
	\draw [thick] (9) -- (12) -- (13);
	\draw [thick] (11) -- (12) -- (8);
	\node at (1.5,-1.5) {$\Phi^+(D_4)$};
\end{tikzpicture}  \vline \; \vline \hspace{-15pt} \begin{tikzpicture}[scale=0.7]
	\node[shape=circle,fill=black,inner sep=1.5] (1) at (0,0) {};
	\node[shape=circle,fill=black,inner sep=1.5] (2) at (1,0) {};
	\node[shape=circle,fill=black,inner sep=1.5] (3) at (0.5,0.5) {};
	\node[shape=circle,fill=black,inner sep=1.5] (4) at (0,1) {};
	\node[shape=circle,fill=black,inner sep=1.5] (5) at (-0.5,1.5) {};
	\node[shape=circle,fill=black,inner sep=1.5] (6) at (-1,2) {};
	\draw [thick] (1) -- (3) -- (2);
	\draw [thick] (3) -- (4) -- (5) -- (6);
	\rotatebox{-45}{\draw [decorate,decoration={brace,amplitude=10pt},yshift=8.5pt] (-2.2,0.5) -- (0.8,0.5) node [black,midway,yshift=15pt]  {\footnotesize $(\ell-1)$};}
	\node at (0,-1.5) {$\Phi^+(I_2(\ell))$};
\end{tikzpicture} \vline \; \; \begin{tikzpicture}[scale=0.5]
	\node[shape=circle,fill=black,inner sep=1.5] (1) at (0,0) {};
	\node[shape=circle,fill=black,inner sep=1.5] (2) at (1,0) {};
	\node[shape=circle,fill=black,inner sep=1.5] (3) at (2,0) {};
	\node[shape=circle,fill=black,inner sep=1.5] (4) at (0.5,0.5) {};
	\node[shape=circle,fill=black,inner sep=1.5] (5) at (1.5,0.5) {};
	\node[shape=circle,fill=black,inner sep=1.5] (6) at (1,1) {};
	\node[shape=circle,fill=black,inner sep=1.5] (7) at (2,1) {};
	\node[shape=circle,fill=black,inner sep=1.5] (8) at (1.5,1.5) {};
	\node[shape=circle,fill=black,inner sep=1.5] (9) at (2.5,1.5) {};
	\node[shape=circle,fill=black,inner sep=1.5] (10) at (1,2) {};
	\node[shape=circle,fill=black,inner sep=1.5] (11) at (2,2) {};
	\node[shape=circle,fill=black,inner sep=1.5] (12) at (1.5,2.5) {};
	\node[shape=circle,fill=black,inner sep=1.5] (13) at (1,3) {};
	\node[shape=circle,fill=black,inner sep=1.5] (14) at (0.5,3.5) {};
	\node[shape=circle,fill=black,inner sep=1.5] (15) at (0,4) {};
	\draw [thick] (1) -- (4) -- (2) -- (5) -- (3);
	\draw [thick] (4) -- (6) -- (5) -- (7) -- (8) -- (6);
	\draw [thick] (7) -- (9) -- (11) -- (8) -- (10) -- (12) -- (11);
	\draw [thick] (12) -- (13) -- (14) -- (15);
	\node at (1,-1.5) {$\Phi^+(H_3)$};
\end{tikzpicture}
\end{center}
\caption{Left: some crystallographic root posets. Right: the non-crystallographic root posets of coincidental type.} \label{fig:root_posets}
\end{figure}
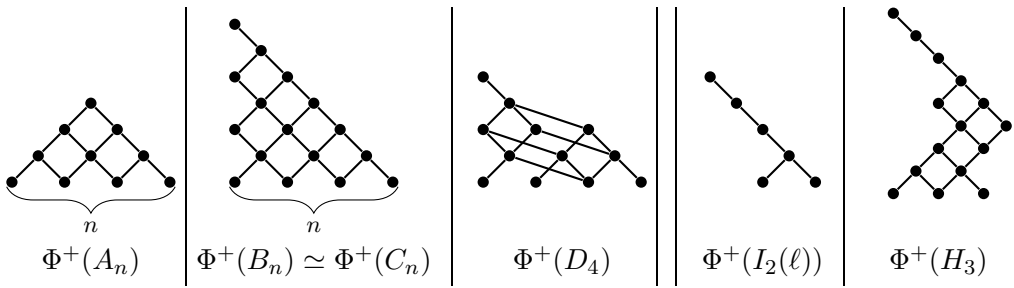

Let $\Phi$ be an irreducible crystallographic root system of rank $n$, and $W$ its Weyl group. We use $\Phi^+$ to denote the positive roots of $\Phi$. We view $\Phi^+$ as a poset, called the \dfn{root poset}, where the partial order is $\alpha \leq \beta$  if and only if $\beta-\alpha$ is a nonnegative sum of simple roots. The poset $\Phi^+$ has $n$ minimal elements (the simple roots) and a unique maximal root (the highest root). It is a graded poset (with the rank function given by height). And it contains important numerical information about $W$. For instance, if $\Phi^+_0,\Phi^+_1,\ldots,\Phi^+_{r(\Phi^+)}$ are the ranks of this poset, then $(\#\Phi^+_0,\#\Phi^+_1,\ldots,\Phi^+_{r(\Phi^+)})$ is a partition and its conjugate partition is $(d_n-1,d_{n-1}-1,\ldots,d_1-1)$, where $d_1\leq d_2 \leq \cdots\leq d_n$ are the degrees of $W$. In particular $r(P)+2=d_n=h$ is the Coxeter number of $W$. Some of these crystallographic root posets are depicted on the left in Figure~\ref{fig:root_posets}.

Now let $\Phi$ be an irreducible non-crystallographic root system. One could naively apply the same definition of partial order to $\Phi^+$, but it would fail to have the desirable features discussed in the last paragraph. Armstrong~\cite{armstrong2009generalized} gave an \emph{ad hoc} construction of root posets $\Phi^+(I_2(\ell))$ and $\Phi^+(H_3)$ which do have these desirable features; these are depicted on the right in Figure~\ref{fig:root_posets}. For $H_4$, the other non-crystallographic root system, there are either many or no analogous root posets, depending on exactly which properties one chooses~\cite{cuntz2015root}.

Among all the complex reflection groups, a special sub-class are the so-called ``coincidental types,'' which are those whose degrees form an arithmetic progression. For the finite Coxeter groups, these are the types $A_n$, $B_n\simeq C_n$, $I_2(\ell)$, and $H_3$. Their corresponding root posets are the \dfn{root posets of coincidental type}. These are all depicted in Figure~\ref{fig:root_posets}. Observe that $\Phi^+(A_n)\simeq S(n)^*$ and $\Phi^+(B_n)\simeq T(n,n)^*$.

For $\Phi^+$ a crystallographic root poset, there is a canonical involutive poset automorphism $\delta\colon \Phi^+ \to \Phi^+$ defined by $\delta(\alpha)  \coloneqq -w_0(\alpha)$ where $w_0\in W$ is the longest element of the Weyl group. For the cases that concern us: $\delta$ is the reflection across the vertical axis of symmetry for $\Phi^+(A_n)$; and $\delta$ is the identity for $\Phi^+(B_n)$. By convention, we define a poset automorphism $\delta\colon \Phi^+ \to \Phi^+$ for $\Phi^+$ a non-crystallographic root poset of coincidental type by: $\delta$ swaps the minimal elements of $\Phi^+(I_2(\ell))$ if $\ell$ is odd, and is the identity if~$\ell$ is even; and $\delta$ is the identity for $\Phi^+(H_3)$.

\begin{figure}
\begin{center}
\begin{tikzpicture}[scale=0.35]
	\node[shape=circle,fill=black,inner sep=1.5] (1) at (0,0) {};
	\node[shape=circle,fill=black,inner sep=1.5] (2) at (-1,1) {};
	\node[shape=circle,fill=black,inner sep=1.5] (3) at (-2,2) {};
	\node[shape=circle,fill=black,inner sep=1.5] (4) at (-3,3) {};
	\node[shape=circle,fill=black,inner sep=1.5] (5) at (-1,3) {};
	\node[shape=circle,fill=black,inner sep=1.5] (6) at (-2,4) {};
	\node[shape=circle,fill=black,inner sep=1.5] (7) at (-3,5) {};
	\node[shape=circle,fill=black,inner sep=1.5] (8) at (-4,6) {};
	\draw [thick] (1) -- (2) -- (3) -- (4) -- (6) -- (7) -- (8);
	\draw [thick] (3) -- (5) -- (6);
	\draw [decorate,decoration={brace,amplitude=10pt},yshift=10pt] (-4,6) -- (-0.9,2.9) node [black,midway,yshift=12pt,xshift=12pt]  {\rotatebox{-45}{\footnotesize $n$}};
	\rotatebox{-45}{\draw [decorate,decoration={brace,mirror,amplitude=10pt},yshift=-10pt] (-4.2,0) -- (0,0) node [black,midway,yshift=-15pt]  {\footnotesize $n$};}
	\node at (-2,-4) {$D(n)$};
\end{tikzpicture} \qquad \vline \qquad
\begin{tikzpicture}[scale=0.3]
	\node[shape=circle,fill=black,inner sep=1.5] (-4) at (3,-3) {};
	\node[shape=circle,fill=black,inner sep=1.5] (-3) at (2,-2) {};
	\node[shape=circle,fill=black,inner sep=1.5] (-2) at (1,-1) {};
	\node[shape=circle,fill=black,inner sep=1.5] (-1) at (2,0) {};
	\node[shape=circle,fill=black,inner sep=1.5] (0) at (-1,1) {};
	\node[shape=circle,fill=black,inner sep=1.5] (1) at (0,0) {};
	\node[shape=circle,fill=black,inner sep=1.5] (2) at (1,1) {};
	\node[shape=circle,fill=black,inner sep=1.5] (3) at (2,2) {};
	\node[shape=circle,fill=black,inner sep=1.5] (4) at (3,3) {};
	\node[shape=circle,fill=black,inner sep=1.5] (6) at (0,2) {};
	\node[shape=circle,fill=black,inner sep=1.5] (7) at (1,3) {};
	\node[shape=circle,fill=black,inner sep=1.5] (8) at (2,4) {};
	\node[shape=circle,fill=black,inner sep=1.5] (11) at (0,4) {};
	\node[shape=circle,fill=black,inner sep=1.5] (12) at (1,5) {};
	\node[shape=circle,fill=black,inner sep=1.5] (13) at (0,6) {};
	\node[shape=circle,fill=black,inner sep=1.5] (14) at (-1,7) {};
	\draw [thick] (-4) -- (-3) -- (-2) -- (-1) -- (2);
	\draw [thick] (2) -- (6) -- (0) -- (1) -- (-2);
	\draw [thick] (1) -- (2) -- (3) -- (4) -- (8) -- (12) -- (13) -- (14);
	\draw [thick] (3) -- (7) -- (6);
	\draw [thick] (7) -- (8);
	\draw [thick] (7) -- (11) -- (12);
	\node at (0,-5) {$\Lambda_{E_6}$};
\end{tikzpicture} \qquad  \vline \qquad \begin{tikzpicture}[scale=0.4]
	\node[shape=circle,fill=black,inner sep=1.5] (1) at (0,0) {};
	\node[shape=circle,fill=black,inner sep=1.5] (2) at (1,0) {};
	\node[shape=circle,fill=black,inner sep=1.5] (3) at (2,0) {};
	\node[shape=circle,fill=black,inner sep=1.5] (4) at (0.5,0.5) {};
	\node[shape=circle,fill=black,inner sep=1.5] (5) at (1.5,0.5) {};
	\node[shape=circle,fill=black,inner sep=1.5] (6) at (1,1) {};
	\node[shape=circle,fill=black,inner sep=1.5] (7) at (2,1) {};
	\node[shape=circle,fill=black,inner sep=1.5] (8) at (1.5,1.5) {};
	\node[shape=circle,fill=black,inner sep=1.5] (9) at (2.5,1.5) {};
	\node[shape=circle,fill=black,inner sep=1.5] (10) at (1,2) {};
	\node[shape=circle,fill=black,inner sep=1.5] (11) at (2,2) {};
	\node[shape=circle,fill=black,inner sep=1.5] (12) at (1.5,2.5) {};
	\node[shape=circle,fill=black,inner sep=1.5] (13) at (1,3) {};
	\node[shape=circle,fill=black,inner sep=1.5] (14) at (0.5,3.5) {};
	\node[shape=circle,fill=black,inner sep=1.5] (15) at (0,4) {};
	\node[shape=circle,fill=black,inner sep=1.5] (-4) at (0.5,-0.5) {};
	\node[shape=circle,fill=black,inner sep=1.5] (-5) at (1.5,-0.5) {};
	\node[shape=circle,fill=black,inner sep=1.5] (-6) at (1,-1) {};
	\node[shape=circle,fill=black,inner sep=1.5] (-7) at (2,-1) {};
	\node[shape=circle,fill=black,inner sep=1.5] (-8) at (1.5,-1.5) {};
	\node[shape=circle,fill=black,inner sep=1.5] (-9) at (2.5,-1.5) {};
	\node[shape=circle,fill=black,inner sep=1.5] (-10) at (1,-2) {};
	\node[shape=circle,fill=black,inner sep=1.5] (-11) at (2,-2) {};
	\node[shape=circle,fill=black,inner sep=1.5] (-12) at (1.5,-2.5) {};
	\node[shape=circle,fill=black,inner sep=1.5] (-13) at (1,-3) {};
	\node[shape=circle,fill=black,inner sep=1.5] (-14) at (0.5,-3.5) {};
	\node[shape=circle,fill=black,inner sep=1.5] (-15) at (0,-4) {};
	\draw [thick] (1) -- (4) -- (2) -- (5) -- (3);
	\draw [thick] (4) -- (6) -- (5) -- (7) -- (8) -- (6);
	\draw [thick] (7) -- (9) -- (11) -- (8) -- (10) -- (12) -- (11);
	\draw [thick] (12) -- (13) -- (14) -- (15);
	\draw [thick] (1) -- (-4) -- (2) -- (-5) -- (3);
	\draw [thick] (-4) -- (-6) -- (-5) -- (-7) -- (-8) -- (-6);
	\draw [thick] (-7) -- (-9) -- (-11) -- (-8) -- (-10) -- (-12) -- (-11);
	\draw [thick] (-12) -- (-13) -- (-14) -- (-15);
	\node at (1.5,-4.75) {$\Lambda_{E_7}$};
\end{tikzpicture}  \qquad  \vline \; \vline \qquad \begin{tikzpicture}[scale=0.65]
	\node[shape=circle,fill=black,inner sep=1.5] (A1) at (0,-1) {};
	\node[shape=circle,fill=black,inner sep=1.5] (B1) at (-1,0) {};
	\node[shape=circle,fill=black,inner sep=1.5] (C1) at (1,0) {};
	\node[shape=circle,fill=black,inner sep=1.5] (A2) at (0,0) {};
	\node[shape=circle,fill=black,inner sep=1.5] (B2) at (-1,1) {};
	\node[shape=circle,fill=black,inner sep=1.5] (C2) at (1,1) {};
	\node[shape=circle,fill=black,inner sep=1.5] (A3) at (0,1) {};
	\node[shape=circle,fill=black,inner sep=1.5] (B3) at (-1,2) {};
	\node[shape=circle,fill=black,inner sep=1.5] (C3) at (1,2) {};
	\draw[thick] (B1)--(A1);
	\draw[thick] (C1)--(A1);
	\draw[thick] (A1)--(A2);
	\draw[thick] (B1)--(B2);
	\draw[thick] (C1)--(C2);
	\draw[thick] (B2)--(A2);
	\draw[thick] (C2)--(A2);
	\draw[thick] (A2)--(A3);
	\draw[thick] (B2)--(B3);
	\draw[thick] (C2)--(C3);
	\draw[thick] (B3)--(A3);
	\draw[thick] (C3)--(A3);
	\node at (0,-2) {$V(n)$};
	\draw [decorate,decoration={brace,amplitude=10pt}] (1.25,2.25) -- (1.25,-0.25) node [black,midway,xshift=15pt]  {$n$};
\end{tikzpicture}
\end{center}
\caption{Left: the other minuscule posets besides the rectangle and shifted staircase. Right: the ``chain of V's.''} \label{fig:minuscule_posets}
\end{figure}
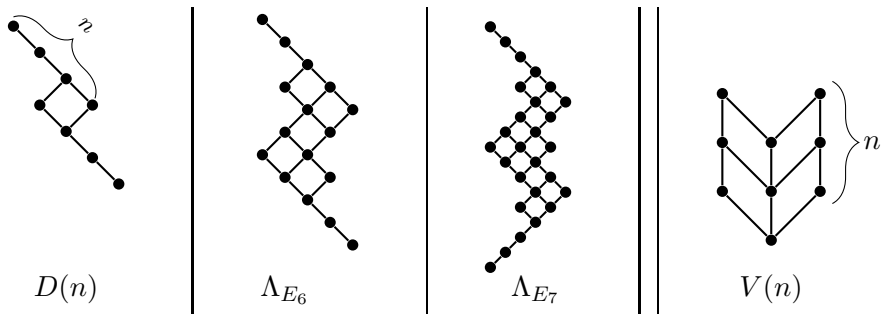

\subsection{Minuscule posets}

The minuscule posets are another family of posets coming from Lie theory, with many remarkable properties. Again, we give only a cursory account; see~\cite[\S4]{hamaker2018doppelgangers} for a detailed treatment.

Let $\mathfrak{g}$ be a simple Lie algebra over $\mathbb{C}$, with $\Phi$ its root system and $W$ its Weyl group. A non-zero (integral, dominant) weight $\lambda$ of $\mathfrak{g}$ is said to be minuscule if the Weyl group acts transitively on the weights of the corresponding highest weight irreducible representation $V^{\lambda}$. In this case, the Weyl orbit $W\lambda$ is a distributive lattice, where the partial order is again given by root order: i.e., $\nu \leq \mu$ means $\mu-\nu$ is a nonnegative sum of simple roots. The \dfn{minuscule poset} corresponding to the minuscule weight $\lambda$ is the poset of join irreducible elements of the distributive lattice $W\lambda$.

If $\lambda$ is minuscule, then it must be equal to some fundamental weight $\omega_i$, and we can also describe the corresponding minuscule poset $P$ as the order filter in~$\Phi^+$ generated by the corresponding simple root $\alpha_i$: i.e., $P=\{\alpha\in \Phi^+\colon \alpha\geq \alpha_i\}$. The minuscule poset $P$ is always graded, and has a unique minimal and a unique maximal element. Furthermore, it has a canonical involutive poset anti-automorphism $\iota\colon P \to P^{*}$ which is induced from the action of multiplication by $w_0$ on $W\lambda$.

The minuscule posets, up to isomorphism, have been classified: they are the \dfn{rectangle} $R(a,b)$, the \dfn{shifted staircase} $DS(n,0)$, the \dfn{``propeller poset''} $D(n)$, and two \dfn{exceptional posets} $\Lambda_{E_6}, \Lambda_{E_7}$ coming from the types $E_6$ and $E_7$. The last three of these are depicted in Figure~\ref{fig:minuscule_posets}.

\subsection{The ``chain of V's''}

The final family of posets with an order polynomial product formula is the \dfn{``chain of $V$'s''}: $V(n)  \coloneqq \begin{tikzpicture}[scale=0.3] \node[shape=circle,fill=black,inner sep=1.5] (B) at (-1,0) {}; \node[shape=circle,fill=black,inner sep=1.5] (C) at (1,0) {}; \node[shape=circle,fill=black,inner sep=1.5] (A) at (0,-1) {}; \draw (B)--(A); \draw (C)--(A); \end{tikzpicture} \times [n]$, the Cartesian product of the $3$-element ``V''-shaped poset $\begin{tikzpicture}[scale=0.3] \node[shape=circle,fill=black,inner sep=1.5] (B) at (-1,0) {}; \node[shape=circle,fill=black,inner sep=1.5] (C) at (1,0) {}; \node[shape=circle,fill=black,inner sep=1.5] (A) at (0,-1) {}; \draw (B)--(A); \draw (C)--(A); \end{tikzpicture}$ and the $n$-element chain $[n]$ (see the right of Figure~\ref{fig:minuscule_posets}). It was first studied by Kreweras and Niederhausen~\cite{kreweras1981solution}. It has a rather different structure than the other examples: for instance, in the other examples each element covers and is covered by at most two elements, but this is not true for $V(n)$. Let us use $\delta\colon V(n)\to V(n)$ to denote the involutive poset automorphism which is reflection across the vertical axis of symmetry of the ``V.''

\section{Order polynomial product formulas} \label{sec:formulas}

In this section we review the order polynomial product formulas for the posets introduced in Section~\ref{sec:posets}, and briefly explain where these formulas come from.

\subsection{Symmetry classes of plane partitions} \label{subsec:sym_classes}

The origin of all these kind of product formulas is MacMahon's investigation of plane partitions, and the subsequent investigation of plane partitions symmetry classes. See~\cite{krattenthaler2016plane} for a complete history.

An $a\times b$ plane partition of height $m$ is an $a\times b$ array $\pi=(\pi_{i,j})_{\substack{1\leq i \leq a,\\ 1\leq j \leq b}}$ of nonnegative integers which is weakly decreasing in rows and columns and for which the largest entry is at most $m$. We denote the set of such plane partitions by $\mathcal{PP}^{m}(a \times b)$. Observe that $\mathcal{PP}^{m}(a \times b)$ is exactly the same as $\mathcal{PP}^{m}(R(a,b))$.

We define the size of a plane partition $\pi \in \mathcal{PP}^{m}(a \times b)$ by $|\pi| \coloneqq \sum_{\substack{1\leq i \leq a, \\ 1 \leq j \leq b}} \pi_{i,j}$. MacMahon~\cite{macmahon1915combinatory} obtained the following celebrated product formula for the size generating function of plane partitions: 
\[\sum_{\pi\in \mathcal{PP}^{m}(a \times b)} q^{|\pi|} = \prod_{i=1}^{a} \prod_{j=1}^{b} \frac{(1-q^{m+i+j-1})}{(1-q^{i+j-1})} \]
This is the $q$-analog $\Omega_P(m;q)$ from Section~\ref{subsec:heuristic} for $P=R(a,b)$.

The order polynomials for other posets beyond $R(a,b)$ arise when considering symmetries of plane partitions. The symmetries relevant to us are transposition $\mathrm{Tr}\colon \mathcal{PP}^{m}(n\times n) \to \mathcal{PP}^{m}(n\times n)$ given by $\mathrm{Tr}(\pi)_{i,j}  \coloneqq \pi_{j,i}$, and complementation $\mathrm{Co}\colon \mathcal{PP}^{m}(a\times b) \to \mathcal{PP}^{m}(a\times b)$ given by $\mathrm{Co}(\pi)_{i,j} \coloneqq m-\pi_{a+1-i,b+1-j}$.

Plane partitions $\pi \in \mathcal{PP}^{m}(n \times n)$ with $\mathrm{Tr}(\pi)=\pi$ are called symmetric; they are evidently in bijection with $\mathcal{PP}^{m}(DS(n,0))$. MacMahon~\cite{macmahon1899partitions} conjectured, and Andrews~\cite{andrews1978plane} and Macdonald~\cite{macdonald1979symmetric} proved, the following product formula for the size generating function of symmetric plane partitions:
\[ \sum_{\pi\in \mathcal{PP}^{m}(n \times n),\mathrm{Tr}(\pi)=\pi} q^{|\pi|} =\prod_{1\leq i < j \leq n}\frac{(1-q^{2(i+j+m-1)})}{(1-q^{2(i+j-1)})} \cdot \prod_{i=1}^{n} \frac{(1-q^{2i+m-1})}{(1-q^{2i-1})}.\]
There is a second $q$-analog for symmetric plane partitions as well. Namely, for a plane partition $\pi \in \mathcal{PP}^{m}(n \times n)$ define $|\pi|'  \coloneqq \sum_{1\leq i \leq j \leq n} \pi_{i,j}$. Then, Bender and Knuth~\cite{bender1972enumeration} conjectured, and Gordon~\cite{gordon1983bender}, Andrews~\cite{andrews1977plane}, and Macdonald~\cite{macdonald1979symmetric} proved:
\[ \sum_{\pi\in \mathcal{PP}^{m}(n \times n),\mathrm{Tr}(\pi)=\pi} q^{|\pi|'} = \prod_{1\leq i \leq j \leq n}\frac{(1-q^{i+j+m-1})}{(1-q^{i+j-1})}.\]
It is this second $q$-analog which is the $\Omega_P(m;q)$ from Section~\ref{subsec:heuristic} for $P=DS(n,0)$.

Plane partitions $\pi \in \mathcal{PP}^{2m}(n \times n)$ with $\mathrm{Tr}(\pi)=\mathrm{Co}(\pi)$ are called transpose-complementary; they are in bijection with $\mathcal{PP}^{m}(S(n))$. These were first enumerated by Proctor (see Theorem~\ref{thm:ap_prod_form}).

Plane partitions $\pi \in \mathcal{PP}^{2m}(n \times n)$ with $\mathrm{Tr}(\pi)=\pi$ and $\mathrm{Co}(\pi)=\pi$ are called symmetric self-complementay; they are in bijection with $\mathcal{PP}^{m}(T(\lfloor n/2 \rfloor, \lceil n/2\rceil) )$. These were again first enumerated by Proctor (see Theorem~\ref{thm:doppelganger_prod_form}).

\subsection{Minuscule posets}

Let $P$ be a poset. For a $P$-partition $\pi \in \mathcal{PP}^m(P)$, define its size to be $|\pi|  \coloneqq \sum_{p \in P}\pi(p)$. Let $F_P(m;q)  \coloneqq \sum_{\pi \in \mathcal{PP}^m(P)} q^{|\pi|}$ denote the size generating function of these $P$-partitions. The basic theory of $P$-partitions (see, e.g., \cite[\S3.15.2]{stanley2012ec1}) says that
\[(1-q)(1-q^2)\cdots(1-q^{\#P}) \cdot \lim_{m\to \infty} F_{P^*}(m;q) = \sum_{L \in \mathcal{L}(P)}q^{\mathrm{maj}(L)},\]
the major index generating function of linear extensions of $P$ (with respect to any fixed natural labeling). We remark that every $d$-complete poset $P$ has a product formula for $\lim_{m\to \infty} F_{P^*}(m;q)$~\cite{proctor2014dcomplete, kim2019hook, naruse2019skew}.

\begin{theorem}[{Proctor~\cite{proctor1984bruhat}}] \label{thm:min_prod_form}
Let $P$ be a minuscule poset. Then
\[ F_P(m;q) =  F_{P^*}(m;q) = \prod_{p \in P}\frac{(1-q^{m+r(p)+1})}{(1-q^{r(p)+1})}, \]
where $r\colon P\to \mathbb{N}$ is the rank function of $P$.
\end{theorem}

Theorem~\ref{thm:min_prod_form} gives the $q$-analog $\Omega_P(m;q)$ from Section~\ref{subsec:heuristic} for $P$ a minuscule poset. Theorem~\ref{thm:min_prod_form} is due to Proctor~\cite{proctor1984bruhat} (although observe that for $R(a,b)$ and $DS(n,0)$ it is equivalent to results just mentioned in Section~\ref{subsec:sym_classes}). To prove this theorem he used Standard Monomial Theory, which explains that $\mathcal{PP}^{m}(P)$ indexes a basis of the representation $V^{m\lambda}$ when $P$ is the minuscule poset corresponding to the minuscule weight $\lambda$. Proctor moreover conjectured that the minuscule posets are the \emph{only} posets which have a product formula for $F_P(m;q)$ of this form.

\subsection{Root posets}

Let $W$ be a finite Coxeter group of rank $n$. Define the $q$-$W$-Catalan number by
\[\mathrm{Cat}(W;q)  \coloneqq \prod_{i=1}^{n} \frac{(1-q^{h+d_i})}{(1-q^{d_i})}\]
where $d_1,\ldots,d_n$ are the degrees of $W$ and $h$ its Coxeter number. It is uniformly known that $\mathrm{Cat}(W; q)$ is a polynomial with nonnegative integer coefficients (see~\cite{bessis2011cyclic}). 
 
An \dfn{order ideal} of a poset $P$ is a downwards-closed subset $I\subseteq P$ (i.e., a subset for which $y\in I$ and $x\leq y\in P$ implies $x\in I$). We use $\mathcal{J}(P)$ to denote the order ideals of $P$. We have a natural identification $\mathcal{J}(P)\simeq \mathcal{PP}^{1}(P)$ where an order ideal corresponds to the indicator function of its complement.

\begin{theorem}[{Cellini--Papi~\cite{cellini2002adnilpotent}, Haiman~\cite{haiman1994conjectures}}] \label{thm:root_j_prod_form}
Let $\Phi$ be a crystallographic root system and $W$ its Weyl group. Then $\#\mathcal{J}(\Phi^+) = \mathrm{Cat}(W;q \coloneqq 1)$.
\end{theorem}

The proof of Theorem~\ref{thm:root_j_prod_form} is uniform; however, there is no known statistic for order ideals of which $\mathrm{Cat}(W;q)$ is the generating function.

Now let $W$ be a finite Coxeter group of coincidental type. Define the $q$-$W$-multi-Catalan number by
\[\mathrm{Cat}(W,m;q)  \coloneqq \prod_{j=0}^{m-1} \prod_{i=1}^{n} \frac{(1-q^{h+d_i+2j})}{(1-q^{d_i+2j})}.\]
It is known, in a case-by-case fashion, that $\mathrm{Cat}(W,m;q)$ is a polynomial in $q$ with nonnegative integer coefficients (for instance, for Type A this follows from~\cite[Theorem~1, Case~`CGI']{proctor1990new}; for other types it can be deduced from consideration of the poset's minuscule doppelg\"{a}nger in the sense of Section~\ref{subsec:doppelganger}).

The multi-Catalan numbers are {\bf not} the same as the more well-known Fuss--Catalan numbers. They first appeared, with this name, in the paper of Ceballos--Labb\'{e}--Stump~\cite{ceballos2014subword} which studied multi-triangulations and the multi-cluster complex. Our interest in these numbers is, however, the following:

\begin{theorem}[{\cite{proctor1983trapezoid, proctor1990new, williams2013cataland}}] \label{thm:root_prod_form}
Let $\Phi$ be a coincidental type root system and~$W$ its corresponding Coxeter group. Then $\Omega_{\Phi^+}(m) = \mathrm{Cat}(W,m;q \coloneqq 1)$.
\end{theorem}

Theorem~\ref{thm:root_prod_form} gives the $q$-analog $\Omega_P(m;q)$ from Section~\ref{subsec:heuristic} for $P$ a root poset of coincidental type. The proof of Theorem~\ref{thm:root_prod_form} is case-by-case, with the difficult cases of $\Phi^+(A_n)$ and $\Phi^+(B_n)$ proved by Proctor~\cite{proctor1990new, proctor1983trapezoid} using representations of the symplectic group. The cases $\Phi^+(I_2(\ell))$ and $\Phi^+(H_3)$ were checked by Williams~\cite{williams2013cataland}. 

\subsection{Doppelg\"{a}ngers} \label{subsec:doppelganger}

Following~\cite{hamaker2018doppelgangers}, we call a pair of posets 
\[(P,Q) \in \{(R(a,b),T(a,b)^*),(DS(5,0),\Phi^+(H_3)),(D(\ell),\Phi^+(I_2(2\ell)))\}\]
a \dfn{minuscule doppelg\"{a}nger pair}.

\begin{theorem}[\cite{proctor1983trapezoid}] \label{thm:doppelganger_prod_form}
Let $(P,Q)$ be a minuscule doppelg\"{a}nger pair. Then we have $\Omega_P(m) = \Omega_Q(m)$ for all $m\geq 1$.
\end{theorem}

Since minuscule posets have product formulas for their order polynomials (Theorem~\ref{thm:min_prod_form}), Theorem~\ref{thm:doppelganger_prod_form} says their doppelg\"{a}ngers do too. The difficult case of Theorem~\ref{thm:doppelganger_prod_form}, the rectangle/trapezoid pair, was yet again established by Proctor~\cite{proctor1983trapezoid} using representations of the symplectic group. The other cases are an easy check. 

In~\cite{hamaker2018doppelgangers}, the authors gave a uniform bijection between $\mathcal{PP}^{m}(P)$ and $\mathcal{PP}^{m}(Q)$ (and also a uniform bijection between $\mathcal{L}(P)$ and $\mathcal{L}(Q)$) for any minuscule doppelg\"{a}nger pair $(P,Q)$. In~\cite{hopkins2019minuscule} it was suggested that minuscule doppelg\"{a}nger pairs are ``very similar.'' We will see more similarities below.

\subsection{Other examples}

\begin{theorem}[{Proctor}] \label{thm:ap_prod_form}
For $P=AP(M,d,\ell)$ an arithmetic progression,
\[ \Omega_{P}(m) = \prod_{\substack{(i,j) \in P, \\ \ell+c(i,j)\leq M-id} } \frac{m+\ell+c(i,j)}{\ell+c(i,j)} \cdot  \prod_{\substack{(i,j) \in P, \\ \ell+c(i,j)> M-id} } \frac{(d+1)m+\ell+c(i,j)}{\ell+c(i,j)}, \]
where $c(i,j)  \coloneqq j-i$ is the content of the box $(i,j)\in P$.
\end{theorem}

Theorem~\ref{thm:ap_prod_form} is due to Proctor. The case $d=1$ has an interpretation in terms of the representation theory of the symplectic group (see~\cite{proctor1988odd}), and the minuscule case $d=0$ has an interpretation in terms of the representation theory of the general linear group. However, for general $d$, Proctor's proof instead manipulates a determinantal formula due to MacMahon (see~\cite[Exercise 7.101]{stanley1999ec2}).

\begin{theorem}[{Kreweras--Niederhausen~\cite{kreweras1981solution}}] \label{thm:kreweras}
For the ``chain of V's,''
\[\Omega_{V(n)}(m) = \frac{\prod_{i=1}^{n}(m+1+i) \prod_{i=1}^{2n}(2m+i+1)}{(n+1)!(2n+1)!}.\]
\end{theorem}

Theorem~\ref{thm:kreweras} is due to Kreweras and Niederhausen~\cite{kreweras1981solution}. Their proof uses some basic $P$-partition theory together with a lot of clever algebraic manipulation and recurrences. Earlier, Kreweras~\cite{kreweras1965classe} obtained the product formula for~$e(V(n))$, which has an interpretation in terms of a $3$ candidate ballot problem.

\begin{theorem}[{Hopkins--Lai~\cite{hopkins2020plane}; Okada~\cite{okada2020intermediate}}] \label{thm:sds}
For the shifted double staircase,
\[ \Omega_{DS(n,k)}(m) = \prod_{1\leq i \leq j  \leq n} \frac{m+i+j-1}{i+j-1} \cdot \prod_{1\leq i \leq j \leq k} \frac{m+i+j}{i+j}.\]
\end{theorem}

The proof of Theorem~\ref{thm:sds} in~\cite{hopkins2020plane} is based on the theory of lozenge tilings of the triangular lattice and the Kuo condensation recurrence technique. In~\cite{okada2020intermediate}, Okada proves Theorem~\ref{thm:sds} using ``intermediate symplectic group''~\cite{proctor1988odd} characters; in fact, he obtains a $q$-analog $\Omega_{DS(n,k)}(m;q)$.

\section{Promotion of linear extensions}

In this section we survey the posets which have good behavior of promotion.

\subsection{Definitions and basics} \label{sec:pro_basics}

Let $P$ be a poset on $n$ elements. For our purposes it is best to represent a linear extension of $P$ as a list $(p_1,p_2,\ldots,p_n)$ of the elements of $P$, each appearing once, for which $p_i \leq p_j$ implies that~$i\leq j$. We then define for each $i=1,\ldots,n-1$ the \dfn{Bender--Knuth involution} $\tau_i\colon \mathcal{L}(P)\to \mathcal{L}(P)$ by 
\[ \tau_i(p_1,\ldots,p_{n})  \coloneqq \begin{cases} (p_1,\ldots,p_{i-1},p_{i+1},p_{i},p_{i+2},\ldots,p_{n}) &\textrm{if $p_{i}$, $p_{i+1}$ incomparable;} \\ (p_1,\ldots,p_{n}) &\textrm{otherwise}. \end{cases} \]
\dfn{Promotion} $\mathrm{Pro}\colon \mathcal{L}(P) \to \mathcal{L}(P)$ is the following composition of these $\tau_i$:
\[ \mathrm{Pro}  \coloneqq \tau_{n-1} \circ \tau_{n-2} \circ \cdots \circ \tau_1. \]
\dfn{Evacuation} $\mathrm{Evac}\colon \mathcal{L}(P) \to \mathcal{L}(P)$ is the following composition of the $\tau_i$:
\[ \mathrm{Evac}  \coloneqq (\tau_1) \circ (\tau_2 \circ \tau_{1}) \circ \cdots \circ (\tau_{n-2} \circ \cdots \circ \tau_{2} \circ \tau_1) \circ (\tau_{n-1} \circ \cdots \circ \tau_{2} \circ \tau_1).\]
There is a duality $\mathcal{L}(P)\to \mathcal{L}(P^*)$ which sends $L=(p_1,\ldots,p_n)$ to $L^*=(p_n,\ldots,p_1)$. Dual evacuation $\mathrm{Evac}^*\colon \mathcal{L}(P) \to \mathcal{L}(P)$ is defined by $\mathrm{Evac}^*(L)  \coloneqq \mathrm{Evac}(L^*)^*$. 

The following are the basic results concerning promotion and evacuation established by Sch\"{u}tzenberger~\cite{schutzenberger1972promotion}; see also the presentation of Stanley~\cite{stanley2009promotion}:

\begin{prop}[{Sch\"{u}tzenberger~\cite{schutzenberger1972promotion}}] \label{prop:pro_basics}
For any poset $P$,
\begin{itemize}
\item $\mathrm{Evac}$ and $\mathrm{Evac}^{*}$ are both involutions;
\item $\mathrm{Evac} \circ \mathrm{Pro}= \mathrm{Pro}^{-1} \circ \mathrm{Evac}$;
\item $\mathrm{Pro}^{\#P} =  \mathrm{Evac}^{*} \circ \mathrm{Evac}$.
\end{itemize}
\end{prop}

Proposition~\ref{prop:pro_basics} explains why $\mathrm{Pro}^{\#P}$ is the ``right'' power of promotion to look at. As mentioned in Section~\ref{subsec:heuristic}, $\mathrm{Pro}^{\#P}$ is ideally the identity, but we will also see interesting examples where it is a non-identity involutive poset automorphism. Let us remark that if $\mathrm{Pro}^{\#P}$ is a poset automorphism, then it must be an involution. Indeed, suppose $\mathrm{Pro}^{\#P}=\delta$ is an automorphism; then by conjugating $\delta=\mathrm{Evac}^{*}\circ \mathrm{Evac}$ by $\mathrm{Evac}$ we get $\delta=\mathrm{Evac} \circ \mathrm{Evac}^{*}$, since evacuation commutes with any automorphism; in other words, we have $\delta=\delta^{-1}$, as claimed.

\subsection{Models}

We now review models (diagrammatic, algebraic, etc.) for promotion and evacuation for certain families of posets. These models lead to a precise understanding of the order and orbit structure of these operators.

\subsubsection{Rotation of noncrossing matchings and webs}

A noncrossing matching of~$[2n]$ is a partition of $[2n]$ into blocks of size $2$ for which there is no pair of crossing blocks. D.~White observed that promotion of standard Young tableaux of $2 \times n$ rectangular shape corresponds to rotation of noncrossing matchings of $[2n]$ (see for instance~\cite[\S8]{rhoades2010cyclic}).

Webs are a class of planar graphs Kuperberg introduced to study the invariant theory of Lie algebras. Khovanov and Kuperberg~\cite{khovanov1999web} (see also Tymoczko~\cite{tymoczko2012simple}) defined a bijection between standard Young tableaux of $3 \times n$ rectangular shape and a subset of $\mathfrak{sl}_3$-webs. Petersen, Pylyavskyy, and Rhoades \cite{petersen2009promotion} showed that under this bijection, promotion of tableaux corresponds to rotation of webs. Russell~\cite{russell2013explicit} (see also Patrias~\cite{patrias2019promotion}) showed that rotation for a broader class of $\mathfrak{sl}_3$-webs corresponds to promotion of certain \emph{semi}standard tableaux of $3\times n$ rectangular shape (see Section~\ref{subsec:semistandard}). Finally, Hopkins and Rubey~\cite{hopkins2020promotion} showed that linear extensions of~$V(n)$ can be encoded as certain $3$-edge-colored $\mathfrak{sl}_3$-webs for which promotion again corresponds to rotation. 

In these correspondences between linear extensions and diagrams it is also possible to show that evacuation corresponds to reflection across a diameter (see, e.g., \cite{patrias2021evacuation}), so that the full dihedral action is apparent.

\subsubsection{Rotation of reduced words and Edelman--Greene-style bijections} \label{subsec:edleman-greene}

Let $W$ be a finite Coxeter group with $\Phi$ its root system. For a reduced word $s_{\alpha_{i_1}}s_{\alpha_{i_2}}\cdots s_{\alpha_{i_p}}$ of the longest word $w_0 \in W$, we define its (twisted) rotation to be $s_{\alpha_{i_2}}\cdots s_{\alpha_{i_p}}s_{\delta(\alpha_{i_1})}$, which is again a reduced word of $w_0$. 

For $W$ of coincidental type, there is an equivariant bijection between linear extensions of $\Phi^+$ under promotion and reduced words of $w_0$ under rotation. For Type~A, this is due to Edelman and Greene~\cite{edelman1987balanced}. For Type~B, it was established by Haiman~\cite{haiman1989mixed} and Kra\'{s}kiewicz~\cite{kraskiewicz1989reduced}. For $I_2(\ell)$ and $H_3$ it was shown by Williams~\cite{williams2013cataland}. 

We believe that evacuation should correspond to (twisted) reflection of the reduced word under these Edelman--Greene-style bijections, but we know of nowhere in the literature where this is explicitly stated.

\subsubsection{Kazhdan--Lusztig cell representations}

The theory of Kazhdan--Lusztig cells gives a canonical basis for any irreducible representation of the symmetric group~$S_n$. Rhoades~\cite{rhoades2010cyclic} showed that for an irreducible symmetric group representation of rectangular shape, the action of the long cycle $c=(1,2,\ldots,n)$ corresponds to promotion of standard tableaux in the Kazhdan--Lusztig basis. Rhoades's result built on an earlier results of Berenstein--Zelevinsky~\cite{berenstein1996canonical} and Stembridge \cite{stembridge1996canonical}, who showed that for \emph{any} irreducible symmetric group representation, the action of the longest word $w_0=n(n-1)\cdots1$ corresponds (up to sign) to evacuation of standard tableaux in the Kazhdan--Lusztig basis.

\subsubsection{Crystals and cactus group actions} \label{subsec:crystals_pro}

The cactus group action on the tensor product of crystals for simple Lie algebra representations, as defined by Henriques and Kamnitzer~\cite{henriques2006crystals}, gives rise to notions of evacuation and promotion acting on the corresponding highest weight words of weight zero: see~\cite{fontaine2014cyclic, westbury2016invariant, pfannerer2020promotion}. In the case of $V$ being the vector representation of $\mathfrak{sl}_k$, this cactus group promotion for weight zero highest weight words of $V^{\otimes kn}$ corresponds to promotion of standard Young tableaux of $k\times n$ rectangular shape, and similarly for evacuation.

\subsubsection{The Wronski map and monodromy} For a list $f_1(z), f_2(z), \ldots, f_d(z)$ of~$d$ linearly independent polynomials in $\mathbb{C}[z]$ of degree at most $n-1$, their Wronskian is the determinant of the $d\times d$ matrix whose rows are the derivatives $f^{(i)}_j(z)$ of these polynomials for $i=0,\ldots,d-1$. Up to scale, the Wronskian depends only the linear span of the $f_j(z)$. We thus obtain the Wronski map from the Grassmannian $\mathrm{Gr}(d,n)$ to projective space $\mathbb{P}^{n-1}$.

Standard Young tableaux of $d \times (n-d)$ rectangular shape index the fibers of the Wronski map. In~\cite{purbhoo2013wronksians}, Purbhoo showed that a certain monodromy action for the Wronski map corresponds to promotion for these tableaux. Moreover, in~\cite{purbhoo2018marvellous}, Purbhoo showed that by restricting to those preimages which lie in either the orthogonal or Lagrangian Grassmannian (under a certain embedding of these Type~B/C Grassmannians into the usual Type~A Grassmannian), one can similarly obtain the action of promotion on standard tableaux of shifted staircase or staircase shape.

\subsubsection{Evacuation of minuscule posets}

The extension of Sch\"{u}tzenberger's theory of jeu de taquin~\cite{schutzenberger1977correspondance} to cominuscule Schubert calculus due to Thomas and Yong~\cite[Lemma 5.2]{thomas2016cominuscule} implies that evacuation of a minuscule poset has a simple description in terms of its canonical anti-automorphism:

\begin{theorem}(\cite{schutzenberger1977correspondance, thomas2016cominuscule}) \label{thm:min_evac}
For a minuscule poset $P$ and $L\in \mathcal{L}(P)$, we have $\mathrm{Evac}(L) = \iota(L)^*$.
\end{theorem}

\subsubsection{Doppelg\"{a}ngers bijections} 

The minuscule doppelg\"{a}nger pairs have the same orbit structure of promotion:

\begin{theorem}[{Haiman~\cite{haiman1989mixed, haiman1992dual}}] \label{thm:doppel_pro}
For $(P,Q)$ a minuscule doppelg\"{a}nger pair, there is a bijection between $\mathcal{L}(P)$ and $\mathcal{L}(Q)$ which commutes with promotion.
\end{theorem}

The bijection for the difficult case of Theorem~\ref{thm:doppel_pro}, the rectangle/trapezoid pair, is due to Haiman~\cite{haiman1989mixed, haiman1992dual}. The other cases are an easy check.

\subsection{Order}

We now review the posets $P$ for which $\mathrm{Pro}^{\#P}$ can be described.

\begin{theorem}[\cite{schutzenberger1977correspondance, edelman1987balanced, haiman1992dual}] \label{thm:pro_shapes}
For $P=R(a,b), T(a,b), DS(n,k)$, $\mathrm{Pro}^{\#P}$ is the identity. For $P=S(n)$, $\mathrm{Pro}^{\#P}$ is transposition.
\end{theorem}

The case of Theorem~\ref{thm:pro_shapes} for the rectangle $R(a,b)$ follows from Sch\"{u}tzenberger's theory of jeu de taquin~\cite{schutzenberger1977correspondance}. The case of the staircase $S(n)$ is due to Edelman--Greene~\cite{edelman1987balanced}. The cases of $T(a,b)$ and $DS(n,k)$ are due to Haiman~\cite{haiman1992dual}, who developed a method which recaptures the rectangle and staircase cases as well. In fact, Haiman and Kim~\cite{haiman1992characterization} showed that the \emph{only} shapes and shifted shapes for which $\mathrm{Pro}^{\#P}$ is the identity or transposition are those appearing in Theorem~\ref{thm:pro_shapes}.

\begin{theorem} \label{thm:min_pro}
For $P$ a minuscule poset, $\mathrm{Pro}^{\#P}$ is the identity.
\end{theorem}

Theorem~\ref{thm:min_pro} follows, for instance, from Theorem~\ref{thm:min_evac}.

\begin{theorem} \label{thm:root_pro}
For $P=\Phi^+$ a root poset of coincidental type, $\mathrm{Pro}^{\#P}=\delta$.
\end{theorem}

Theorem~\ref{thm:root_pro} follows, for instance, from the reduced word bijections (see the discussion in Section~\ref{subsec:edleman-greene}).

\begin{theorem}[{Hopkins--Rubey~\cite{hopkins2020promotion}}]
For $P=V(n)$, $\mathrm{Pro}^{\#P}=\delta$.
\end{theorem}

\subsection{Orbit structure}

We now discuss the orbit structure of promotion for the posets with good promotion behavior.

The cyclic sieving phenomenon of Reiner--Stanton--White~\cite{reiner2004cyclic} provides a very compact way to record the orbit structure of a cyclic group action. Recall that if~$X$ is a combinatorial set, $C=\langle c \rangle$ is a cyclic group of order $n$ acting on $X$, and $f(q) \in \mathbb{N}[q]$ is a polynomial with nonnegative integer coefficients, then we say that the triple $(X,C,f(q))$ exhibits \dfn{cyclic sieving} if for all integers $k$,
\[ \#\{x\in X\colon c^k(x)=x\}=f(q \coloneqq \zeta^k),\]
where $\zeta  \coloneqq e^{2\pi i/n}$ is a primitive $n$th root of unity. As mentioned in Section~\ref{subsec:heuristic}, cyclic sieving phenomena (CSPs) where the polynomial has a simple product formula are especially valuable, because they imply a product formula for every symmetry class.

For the minuscule posets there is a beautiful such CSP:

\begin{theorem}[\cite{rhoades2010cyclic, purbhoo2018marvellous, sekheri2014CSP}] \label{thm:min_pro_csp}
Let $P$ be a minuscule poset. Then $(\mathcal{L}(P),\langle \mathrm{Pro} \rangle, f(q))$ exhibits cyclic sieving, where
\[ f(q)  \coloneqq \sum_{L \in \mathcal{L}(P)}q^{\mathrm{maj}(L)} = (1-q)(1-q^2)\cdots (1-q^{\#P})\cdot \prod_{p \in P} \frac{1}{(1-q^{r(p)+1})}.\]
\end{theorem}

The case of Theorem~\ref{thm:min_pro_csp} for the rectangle $R(a,b)$ is due to Rhoades~\cite{rhoades2010cyclic}. The exceptional cases $\Lambda_{E_6}, \Lambda_{E_7}$ are a finite check which is easily carried out by computer. The case $D(n)$ is trivial. For the shifted staircase $DS(n,0)$, Purbhoo~\cite[Theorem~5.1(i)]{purbhoo2018marvellous} gave an interpretation of this fixed point count in terms of ribbon tableaux; and in~\cite{sekheri2014CSP} it was verified that this matches the CSP evaluation. It would be preferable to have a uniform proof of Theorem~\ref{thm:min_pro_csp}.

For the other posets with good promotion behavior, we have conjectural CSPs:

\begin{conj} \label{conj:root_pro_csp}
Let $\Phi$ be a rank $n$ root system of coincidental type, $W$ its Coxeter group, and $h$ its Coxeter number. Let $C=\langle c\rangle\simeq \mathbb{Z}/nh\mathbb{Z}$ act on $\mathcal{L}(\Phi^+)$ via $c(L)  \coloneqq  \mathrm{Pro}(L)$. Then $(\mathcal{L}(\Phi^+),C,f(q))$ exhibits cyclic sieving, where
\[f(q)  \coloneqq (1-q^2)(1-q^4)\cdots (1-q^{nh})\cdot \lim_{m\to \infty} \mathrm{Cat}(W,m;q) \in \mathbb{N}[q].\]
\end{conj}

The cases $\Phi=B_n,H_3,I_2(2\ell)$ of Conjecture~\ref{conj:root_pro_csp} follow from Theorems~\ref{thm:doppel_pro} and~\ref{thm:min_pro_csp}. The case $\Phi=I_2(2\ell+1)$ is trivial. Hence, the only open case is $\Phi=A_n$, for which N.~Williams conjectured this CSP in a different but equivalent form c.~2010. Purbhoo~\cite[Theorem~5.1(ii)]{purbhoo2018marvellous} again gave an interpretation of this fixed point count in terms of ribbon tableaux, so possibly this conjecture could be resolved as in~\cite{sekheri2014CSP}.

\begin{conj}[{Hopkins--Rubey~\cite{hopkins2020promotion}}]
For all $n\geq 1$, the rational expression
\[ f(q)  \coloneqq \frac{\prod_{i=1}^{3n} (1-q^{2i})}{\prod_{i=2}^{2n+1}(1-q^i)\prod_{i=2}^{n+1} (1-q^{2i}) }\]
is in $\mathbb{N}[q]$, and $(\mathcal{L}(V(n)),\langle \mathrm{Pro} \rangle, f(q))$ exhibits cyclic sieving.
\end{conj}

\begin{conj}
For all $1 \leq k \leq n$, the rational expression
\[ f(q)  \coloneqq \frac{\prod_{i=1}^{n(n+1)/2 + k(k+1)/2} (1-q^{i})}{\prod_{1\leq i \leq j \leq n}(1-q^{i+j-1})\prod_{1\leq i 
\leq j \leq k} (1-q^{i+j})} \]
is in $\mathbb{N}[q]$, and $(\mathcal{L}(DS(n,k)),\langle \mathrm{Pro} \rangle, f(q))$ exhibits cyclic sieving.
\end{conj} 

Finally, we note that for evacuation acting on the linear extensions $\mathcal{L}(P)$ of \emph{any} poset $P$, there is a CSP using the (co)major index generating function~\cite[\S3]{stanley2009promotion}.

\section{Rowmotion of order ideals and \texorpdfstring{$P$}{P}-partitions}

In this section we survey the posets which have good behavior of rowmotion.

\subsection{Definitions and basics} \label{sec:row_basics}

Let $P$ be a poset. As discussed in Section~\ref{sec:intro}, rowmotion was originally defined as an action on $\mathcal{J}(P)$. But, following~\cite{einstein2013combinatorial}, we will right away define it as a piecewise-linear action on $\mathcal{PP}^{m}(P)$ for any $m$. For $p\in P$ we define the \dfn{(piecewise-linear) toggle} $\tau_p\colon \mathcal{PP}^{m}(P) \to \mathcal{PP}^{m}(P)$ by
\[\tau_p(\pi)(q)  \coloneqq \begin{cases} \pi(q) &\textrm{if $p\neq q$}; \\ \min(\{\pi(r)\colon p \lessdot r\}) + \max(\{\pi(r)\colon r\lessdot p\}) - \pi(p) &\textrm{if $p=q$},\end{cases}\]
where $\min(\varnothing) \coloneqq m$ and $\max(\varnothing) \coloneqq 0$. \dfn{Rowmotion} $\mathrm{Row}\colon \mathcal{PP}^{m}(P)\to \mathcal{PP}^{m}(P)$ is then the following composition of these $\tau_i$:
\[\mathrm{Row}  \coloneqq \tau_{p_1}\circ \tau_{p_2} \circ \cdots \circ \tau_{p_n},\]
where $(p_1,\ldots,p_n)$ is any linear extension of $P$. The traditional case of order ideal rowmotion is recovered via the identification $\mathcal{J}(P)\simeq\mathcal{PP}^{1}(P)$.

Rowmotion only ever has good behavior when $P$ is graded (indeed, the name ``rowmotion'' indicates toggling ``row-by-row,'' i.e., ``rank-by-rank''). So from now on assume $P$ is graded and $P_0, P_1, \ldots, P_{r(P)}$ are its ranks. Define $\tau_i  \coloneqq \prod_{p\in P_i}\tau_p$ for $i=0,\ldots,r(P)$ (these toggles all commute, so this product makes sense). Observe that $\mathrm{Row} = \tau_{0} \circ \cdots \circ \tau_{r(P)-1} \circ \tau_{r(P)}$. In analogy with promotion/evacuation, let us then define \dfn{rowvacuation} $\mathrm{Rvac}\colon \mathcal{PP}^{m}(P)\to \mathcal{PP}^{m}(P)$ by
\[\mathrm{Rvac}  \coloneqq (\tau_{r(P)}) \circ (\tau_{r(P)-1} \circ \tau_{r(P)}) \circ \cdots \circ (\tau_{1} \circ \cdots \circ \tau_{r(P)-1} \circ \tau_{r(P)}) \circ (\tau_{0} \circ \cdots \circ \tau_{r(P)-1} \circ \tau_{r(P)}) \]
There is a duality $\mathcal{PP}^{m}(P) \to \mathcal{PP}^{m}(P^*)$ given by $\pi^{*}(p)  \coloneqq m-\pi(p)$. Dual rowvacuation $\mathrm{Rvac}^{*}\colon \mathcal{PP}^{m}\to \mathcal{PP}^{m}(P)$ is defined by $\mathrm{Rvac}^{*}(\pi)  \coloneqq \mathrm{Rvac}(\pi^*)^*$. 

The arguments from~\cite{stanley2009promotion} which establish Proposition~\ref{prop:pro_basics} are very formal (they only use that the~$\tau_i$ are involutions and that $\tau_i$ and $\tau_j$ commute if $|i-j|\geq 2$) and apply here as well, so that we have:

\begin{prop} \label{prop:row_basics}
For any graded poset $P$,
\begin{itemize}
\item $\mathrm{Rvac}$ and $\mathrm{Rvac}^{*}$ are both involutions;
\item $\mathrm{Rvac} \circ \mathrm{Row}= \mathrm{Row}^{-1} \circ \mathrm{Rvac}$;
\item $\mathrm{Row}^{r(P)+2} =  \mathrm{Rvac}^{*} \circ \mathrm{Rvac}$.
\end{itemize}
\end{prop}

Proposition~\ref{prop:row_basics} explains why $\mathrm{Row}^{r(P)+2}$ is the ``right'' power of rowmotion to look at. Again, ideally $\mathrm{Row}^{r(P)+2}$ is the identity, but in some interesting cases it is an involutive poset automorphism. The same argument as for promotion shows that if $\mathrm{Row}^{r(P)+2}$ is a poset automorphism, it must be an involution.

{\bf N.B.}: often ``left-to-right'' toggling of order ideals or $P$-partitions is studied, as opposed to the ``top-to-bottom'' toggling of rowmotion. This left-to-right toggling is routinely called ``promotion,'' but to avoid confusion with promotion of linear extensions, we will not use that term. In practice, the techniques of~\cite{striker2012promotion} can always be used to show that left-to-right toggling is conjugate to top-to-bottom toggling.

\subsection{Models}

We now review models for rowmotion and rowvacuation for certain families of posets. As with the models for promotion and evacuation, these models give us a precise understanding of the order and orbit structure.

\subsubsection{Rotation of binary words, and parabolic cosets of Weyl groups}

Under the ``Stanley--Thomas word'' bijection (see~\cite{stanley2009promotion, propp2015homomesy}), rowmotion of $\mathcal{J}(R(a,b))$ corresponds to rotation of binary words with $a$ $1$'s and $b$ $0$'s. Extending this description, Rush and Shi~\cite{rush2013orbits} showed that if $P$ is the minuscule poset corresponding to the minuscule weight~$\lambda$, then under the natural isomorphism $\mathcal{J}(P) \simeq W/W_{J}$, where $W_J$ is the parabolic subgroup of the Weyl group $W$ stabilizing $\lambda$, the action of rowmotion is conjugate to the action of a Coxeter element $c\in W$.

\subsubsection{Kreweras complementation for noncrossing partitions, in all types}

A noncrossing partition of $[n]$ is a set partition of $[n]$ for which there is no pair of crossing blocks. The noncrossing partitions form a lattice, and the Kreweras complementation is an operator acting on this lattice which has order $2n$ (its square is rotation). 

Now let $W$ be a finite Coxeter group and $\Phi$ its root system. By analogy, the noncrossing partitions for $W$ are the elements less than some fixed Coxeter element $c\in W$ in absolute order, and Kreweras complementation is the operator $w\mapsto cw^{-1}$ on these noncrossing partitions. Proving conjectures of Panyushev~\cite{panyushev2009orbits} and Bessis--Reiner~\cite{bessis2011cyclic}, Amstrong, Stump, and Thomas~\cite{armstrong2013uniform} showed that when $W$ is a Weyl group, rowmotion of $\mathcal{J}(\Phi^+)$ is in equivariant bijection with Kreweras complementation of the noncrossing partitions of $W$ (and for non-crystallographic types, see~\cite{cuntz2015root}).

In~\cite{defant2021symmetry}, rowvacuation of $\mathcal{J}(\Phi^+)$ is studied. In particular, it is shown there that the bijection of Amstrong--Stump--Thomas between $\mathcal{J}(\Phi^+)$ and the noncrossing partitions of $W$ transports the action of $\mathrm{Row}^{-1} \circ \mathrm{Rvac}$ to the involution $w \mapsto g w^{-1} g^{-1}$ (for the appropriate involution $g \in W$ depending on $c$). For the diagrammatic models of noncrossing partitions, this involution is a reflection across a diameter. In~\cite{hopkins2020birational} it is shown that, in Type~A, rowvacuation of $\mathcal{J}(\Phi^+)$ is the same as the ``Lalanne--Kreweras involution'' on Dyck paths.

\subsubsection{Promotion and evacuation of semistandard Young tableaux} \label{subsec:semistandard}

A semistandard Young tableau (SSYT) of shape $\lambda$ is a filling of the boxes of $\lambda$ with positive integers that is weakly increasing in rows and strictly increasing in columns. Let $\mathrm{SSYT}(\lambda,k)$ denote the set of semistandard tableaux of shape~$\lambda$ with entries at most $k$. For $i=1,\ldots,k-1$, the $i$th Bender--Knuth involution is an operator on $\mathrm{SSYT}(\lambda,k)$ which exchanges the number of $i$'s and $(i+1)$'s in a tableau. Promotion and evacuation can then be defined as operators on $\mathrm{SSYT}(\lambda,k)$ which are the appropriate compositions of these Bender--Knuth involutions.

Promotion and evacuation of semistandard Young tableaux are well-studied actions, with several algebraic guises (e.g., see Section~\ref{sec:ssyt_canonical} below). In particular, promotion of SSYTs of rectangular shape is used in the combinatorial definition of affine Type~A crystals, or more precisely, the so-called ``Kirillov--Reshetikhin crystals;'' see~\cite{shimozono2002affine, bandlow2010uniqueness}. In some sense promotion corresponds to the cyclic symmetry of the affine Dynkin diagram.

It is known that rowmotion acting on $\mathcal{PP}^{m}(R(a,b))$ is in equivariant bijection with promotion acting on $\mathrm{SSYT}(m^a,a+b)$ (see, e.g.,~\cite{kirillov1995groups} and~\cite[Appendix A]{hopkins2019cyclic}). And likewise for rowvacuation and evacuation.

\subsubsection{Canonical bases, from quantum groups and cluster algebras} \label{sec:ssyt_canonical}

The theory of quantum groups gives canonical bases for Lie group representations. Rhoades \cite{rhoades2010cyclic} (see also~\cite{rush2020restriction}) showed that for an irreducible representation of the general linear group of rectangular shape, the action of the long cycle on the dual canonical basis corresponds to the action of promotion on rectangular SSYTs. Again, Rhoades's result built on an earlier results of Berenstein--Zelevinsky~\cite{berenstein1996canonical} and Stembridge~\cite{stembridge1996canonical}, who showed that for \emph{any} irreducible representation of the general linear group, the action of the longest word corresponds (up to sign) to evacuation of tableaux in the dual canonical basis. Thanks to the discussion in Section~\ref{subsec:semistandard} above, this means that these actions are conjugate to romotion and rowvacuation of $\mathcal{PP}^{m}(R(a,b))$.

The theory of cluster algebras also gives canonical bases. Shen and Weng \cite{shen2018cyclic} showed that the action of the cyclic shift on the theta basis of the coordinate ring of the Grassmannian -- which is a canonical basis coming from its structure as a cluster algebra -- is also conjugate to the action of rowmotion on $\mathcal{PP}^{m}(R(a,b))$.

Recently, Gao, Lam, and Xu~\cite{gao2022electrical} introduced the grove algebra, a version of the coordinate ring of the Lagrangian Grassmannian. They conjectured the existence of an electrical canonical basis for this grove algebra, which comes with a cyclic action. We believe that this cyclic action on electrical canonical basis elements should be in equivariant bijection with rowmotion of $\mathcal{PP}^{m}(\Phi^+(A_n))$.

\subsubsection{Crystals and cactus group actions} \label{sec:crystals}

Recall the cactus group promotion of highest weight words of weight zero for tensor products of crystals discussed in Section~\ref{subsec:crystals_pro}. As explained in~\cite[Example 2.4]{pfannerer2020promotion}, for tensor products of the spin representation of the spin group, these weight zero highest weight words correspond to fans of Dyck paths, which are in bijection with $\mathcal{PP}^{m}(\Phi^+(A_n))$. Moreover, using the techniques of~\cite{pfannerer2020promotion} it can be shown that cactus group promotion of these words is in equivariant bijection with rowmotion of $\mathcal{PP}^{m}(\Phi^+(A_n))$, and similarly for evacuation and rowvacuation. Promotion of fans of Dyck paths is studied in detail in the recent paper~\cite{pappe2022promotion}.

\subsubsection{Quiver representations and reflection functors}

In~\cite{garver2018minuscule}, Garver, Patrias, and Thomas studied minuscule posets and their $P$-partitions from the perspective of quiver representations. Fixing a quiver $Q$ whose underlying graph is a Dynkin diagram, and a node $i$ of this Dynkin diagram corresponding to a minuscule weight~$\omega_i$, they showed that the Jordan form of a generic nilpotent endomorphism gives a bijection from representations $X$ of $Q$ with support at $i$ to $P$-partitions for the corresponding minuscule poset~$P$. Moreover, they described the piecewise-linear toggles in terms of reflection functors. In this way, they were able to analyze rowmotion for $P$-partitions of minuscule posets using quiver representations.

\subsubsection{Rowvacuation of minuscule posets}

Rowvacuation of a minuscule poset has a simple description in terms of its canonical anti-automorphism:

\begin{theorem}[{\cite{grinberg2015birational2, okada2020birational}}] \label{thm:min_rvac}
For a minuscule poset $P$ and $\pi \in \mathcal{PP}^{m}(P)$, we have $\mathrm{Rvac}(\pi) = \iota(\pi)^*$.
\end{theorem}  

The case $P=R(a,b)$ of Theorem~\ref{thm:min_rvac} was proved by Grinberg--Roby~\cite{grinberg2015birational2}; the rest of the theorem was proved, in a case-by-case manner, by Okada~\cite{okada2020birational}. Those authors described their results in terms of ``reciprocity'' of rowmotion, but it is easy to translate their results to this statement about rowvacuation.

\subsubsection{Doppelg\"{a}ngers bijections}

We conjectured that minuscule doppelg\"{a}nger pairs have the same orbit structure of rowmotion: 

\begin{conj}[{\cite{hopkins2019minuscule}}] \label{conj:doppel_row}
For $(P,Q)$ a minuscule doppelg\"{a}nger pair, there is a bijection between $\mathcal{PP}^{m}(P)$ and $\mathcal{PP}^{m}(Q)$ which commutes with rowmotion.
\end{conj}

The case $m=1$ of Conjecture~\ref{conj:doppel_row} was proved in~\cite{dao2019trapezoid}, using the bijection of~\cite{hamaker2018doppelgangers}. In the very recent preprint~\cite{johnson2023plane}, Johnson and Liu proved the main case of Conjecture~\ref{conj:doppel_row} involving the rectangle-trapezoid doppelg\"{a}nger pair.

\subsubsection{Symmetry classes of plane partitions}

Grinberg and Roby~\cite{grinberg2015birational2} explained how rowmotion for the three ``triangular'' posets $DS(n,0)$, $\Phi^+(A_n)$, and $\Phi^+(B_n)$ can be understood by imposing symmetries on the rectangle: 

\begin{lemma}[{Grinberg--Roby~\cite{grinberg2015birational2}; see  also~\cite[\S5]{hopkins2019cyclic}}] \label{lem:row_syms} \hfill
\begin{itemize}
\item There is a $\mathrm{Row}$-equivariant bijection between $\mathcal{PP}^{m}(DS(n,0))$ and the subset $\pi \in \mathcal{PP}^{m}(n\times n)$ with $\mathrm{Tr}(\pi)=\pi$.
\item There is a $\mathrm{Row}$-equivariant bijection between $\mathcal{PP}^{m}(\Phi^+(A_n))$ and the subset of $\pi \in \mathcal{PP}^{2m}((n+1)\times (n+1))$ with $\mathrm{Row}^{n+1}(\pi)=\mathrm{Tr}(\pi)$.
\item There is a $\mathrm{Row}$-equivariant bijection between $\mathcal{PP}^{m}(\Phi^+(B_n))$ and the subset of $\pi \in \mathcal{PP}^{2m}(2n \times 2n)$ with $\mathrm{Tr}(\pi)=\pi$ and $\mathrm{Row}^{2n}(\pi)=\pi$.
\end{itemize}
\end{lemma}

We note that a similar ``triangle-into-rectangle'' embedding, but for linear extensions rather than $P$-partitions, was studied by Pon and Wang~\cite{pon2011promotion}.

\subsection{Order}

We now review the posets $P$ for which $\mathrm{Row}^{r(P)+2}$ acting on $\mathcal{PP}^{m}(P)$ can be described.

\begin{theorem}[\cite{grinberg2015birational2,garver2018minuscule,okada2020birational}] \label{thm:min_row}
For $P$ a minuscule poset, we have that $\mathrm{Row}^{r(P)+2}$ is the identity.
\end{theorem}

\begin{theorem}[\cite{grinberg2015birational2}] \label{thm:root_row}
For $P=\Phi^+$ a root poset of coincidental type, we have $\mathrm{Row}^{r(P)+2}=\delta$.
\end{theorem}

Theorems~\ref{thm:min_row} and~\ref{thm:root_row} were essentially proved, in a case-by-case fashion, by Grinberg--Roby~\cite{grinberg2015birational2}; the only case they could not address was $\Lambda_{E_7}$, which was resolved in~\cite{garver2018minuscule} and~\cite{okada2020birational}.

Conjecture~\ref{conj:doppel_row} and Theorem~\ref{thm:min_row} together imply that the trapezoid should have $\mathrm{Row}^{r(P)+2}$ equal to the identity (this was also conjecture by N.~Williams~\cite[Conjecture~75]{grinberg2015birational2}). But beyond the case $m=1$ this remains open.

The only other poset which apparently has good (piecewise-linear) rowmotion behavior is $V(n)$:

\begin{conj} \label{conj:chain_of_vs_rowmotion}
For $P=V(n)$, we have $\mathrm{Row}^{r(P)+2}=\delta$.
\end{conj}

In the very recent preprint~\cite{adenbaum2023order}, Adenbaum has proved Conjecture~\ref{conj:chain_of_vs_rowmotion}.

\subsection{Orbit structure}

There are very nice conjectural CSPs for all the posets with good $P$-partition rowmotion behavior:

\begin{conj}[{\cite{hopkins2019minuscule}}] \label{conj:min_row_csp}
Let $P$ be a minuscule poset. Then the triple $(\mathcal{PP}^{m}(P),\langle \mathrm{Row} \rangle,F_P(m;q))$ exhibits cyclic sieving.
\end{conj}

\begin{conj}[{\cite{hopkins2019minuscule}}] \label{conj:root_row_csp}
Let $\Phi$ be a root system of coincidental type, $W$ its Coxeter group, and $h$ its Coxeter number. Let $C=\langle c\rangle\simeq \mathbb{Z}/2h\mathbb{Z}$ act on $\mathcal{PP}^{m}(\Phi^+)$ via $c(\pi) \coloneqq \mathrm{Row}(\pi)$. Then $(\mathcal{PP}^{m}(\Phi^+),C,\mathrm{Cat}(W,m;q))$ exhibits cyclic sieving.
\end{conj}

The case $m=1$ of Conjecture~\ref{conj:min_row_csp} was proved by Rush--Shi~\cite{rush2013orbits}. The case $m=1$ of Conjecture~\ref{conj:root_row_csp} was proved by Armstrong--Stump--Thomas~\cite{armstrong2013uniform} (and in fact they showed this for \emph{any} root system; see also~\cite{cuntz2015root}). The case of Conjecture~\ref{conj:min_row_csp} for the rectangle $R(a,b)$ was proved by Rhoades~\cite{rhoades2010cyclic} (and later, Shen--Weng~\cite{shen2018cyclic}). All other cases of these conjectures are open. We note the Type A case of Conjecture~\ref{conj:root_row_csp} was conjectured also by J.~Propp c.~2016.

Via Lemma~\ref{lem:row_syms}, the main remaining cases of Conjectures~\ref{conj:min_row_csp} and~\ref{conj:root_row_csp} can be translated into statements about the numbers of fixed points of various subgroups of $\langle \mathrm{Row}, \mathrm{Tr}\rangle$ acting on $\mathcal{PP}^{m}(n\times n)$ being counted by CSP-type evaluations. In~\cite{hopkins2019cyclic}, it was shown, building off the work of Rhoades~\cite{rhoades2010cyclic}, that the number of fixed points of any \emph{element} of $\langle \mathrm{Row}, \mathrm{Tr}\rangle$ acting on $\mathcal{PP}^{m}(n\times n)$ is given by a CSP-type evaluation of polynomial with a nice product formula.

Conjectures~\ref{conj:doppel_row} and~\ref{conj:min_row_csp} together describe the orbit structure of rowmotion for minuscule doppelg\"{a}ngers (and their claims are consistent with Conjecture~\ref{conj:root_row_csp}).

Finally, for $V(n)$ we conjecture:

\begin{conj}
For all $n,m \geq 1$, the rational expression
\[ f(q)  \coloneqq \prod_{i=2}^{2n+1} \frac{(1-q^{2m+i})}{(1-q^i)} \cdot \prod_{i=2}^{n+1}\frac{(1-q^{2m+2i})}{(1-q^{2i})}\] 
is in $\mathbb{N}[q]$, and $(\mathcal{PP}^{m}(V(n)),\langle \mathrm{Pro} \rangle, f(q))$ exhibits cyclic sieving.
\end{conj}

\bibliographystyle{amsalpha}

\bibliography{hopkins}{}

\end{document}